\documentclass[11pt,reqno]{amsart}

\usepackage{lmodern}
\usepackage[colorlinks=true]{hyperref}
\usepackage{amsmath,amssymb,mathrsfs,amsthm, mathtools}
\usepackage[mathscr]{euscript}
\usepackage[utf8]{inputenc}
\usepackage[T1]{fontenc}
\usepackage{mathabx}

\newtheoremstyle{exampstyle}
{10pt} % Space above
{10pt} % Space below
{\it } % Body font
{} % Indent amount
{\bfseries} % Theorem head font
{.} % Punctuation after theorem head
{.5em} % Space after theorem head
{} % Theorem head spec (can be left empty, meaning `normal')

\theoremstyle{exampstyle}
\newtheorem{theorem}{Theorem}[section]

\newcommand{\pv}{\text{p.v.}\int_{\mathbb{S}^1}}

\def\sideremark#1{\ifvmode\leavevmode\fi\vadjust{\vbox to0pt{\vss
\hbox to0pt{\hskip\hsize\hskip1em
\vbox{\hsize3cm\tiny\raggedright\pretolerance10000
\noindent #1\hfill}\hss}\vbox to8pt{\vfil}\vss}}}

\def\bbS{{\mathbb S^1}}

\def\px{\partial_x}
\def\pt{\partial_t}

\usepackage{scalerel,stackengine}
\stackMath
\newcommand\reallywidehat[1]{%
\savestack{\tmpbox}{\stretchto{%
  \scaleto{%
    \scalerel*[\widthof{\ensuremath{#1}}]{\kern-.6pt\bigwedge\kern-.6pt}%
    {\rule[-\textheight/2]{1ex}{\textheight}}%WIDTH-LIMITED BIG WEDGE
  }{\textheight}%
}{0.5ex}}%
\stackon[1pt]{#1}{\tmpbox}%
}

\newcommand{\w}[1]{{\reallywidehat{#1}}}

\def\q{\quad}
\def\qq{\qquad}

\def\bbZ{{\mathbb{Z}}}

\def\pat{\partial_t}
\def\pa{\partial}

\def\pa{\partial}
\def\eps{\varepsilon}

\def\ds{\displaystyle}

\def\intS{\int_{\bbS}}

\newcommand{\m}[1]{\mathcal{#1}}

\newcommand{\pare}[1]{\left( #1 \right)}
\newcommand{\norm}[1]{\left\| #1 \right\|}
\newcommand{\av}[1]{\left| #1 \right|}
\newcommand{\bra}[1]{\left[ #1 \right]}

\renewcommand{\t}[1]{\text{#1}}

\numberwithin{equation}{section}
\usepackage{geometry}
\usepackage{enumitem}
\setlist[itemize,1]{leftmargin=\dimexpr 20pt}

\usepackage{cite} %\def\citepunct{], [} \def\citedash{]--[}

\title{New functional inequalities with applications to the \emph{arctan}-fast diffusion equation}

\author[R. Granero-Belinch\'{o}n]{Rafael Granero-Belinch\'{o}n}

\address{Departamento  de  Matem\'aticas,  Estad\'istica  y  Computaci\'on,  Universidad  de Cantabria.  Avda.  Los  Castros  s/n,  Santander,  Spain.}
\email{rafael.granero@unican.es}

\author[M. Magliocca]{Martina Magliocca}

\address{Departamento  de  Análisis Matemático,  Universidad  de Sevilla.  C/Tarfia s/n, Campus Reina Mercedes
41012 Sevilla, Spain.}
\email{mmagliocca@us.es}

\author[A. Ortega]{Alejandro Ortega}

\address{Departamento de Matem\'aticas Fundamentales, Facultad de Ciencias, UNED, 28040 Madrid, Spain}
\email{alejandro.ortega@mat.uned.es}

\date{\today}

\allowdisplaybreaks
\begin{document}

\begin{abstract}
In this paper, we prove a couple of new nonlinear functional inequalities of Sobolev type akin to the logarithmic Sobolev inequality. In particular, one of the inequalities reads $$\int_{\bbS}\arctan\pare{\frac{\pa_x u}{u}}\pa_xu \,dx\geq \arctan\pare{\norm{u}_{\dot{W}^{1,1}(\bbS)}}\norm{u}_{\dot{W}^{1,1}(\bbS)}.$$ Then, these inequalities are used in the study of the nonlinear \emph{arctan}-fast diffusion equation
$$
\pt u-\px\arctan\left(\frac{\px u}{u}\right)=0.
$$
For this highly nonlinear PDE we establish a number of well-posedness results and qualitative properties.
\end{abstract}

\keywords{Arctan-fast diffusion equation, Global solutions, Functional inequalities}

\subjclass[2010]{35B09, 35K20, 35K55, 47J20}

\maketitle
{\small
\tableofcontents}

\section{Introduction and main results}
Functional inequalities are at the core of functional analysis and partial differential equations. In this paper, we prove the following nonlinear Sobolev inequality
$$
\int_{\bbS}\arctan\pare{\frac{\av{\pa_x u}}{u}}\av{\pa_xu}\,dx\geq \arctan\pare{\norm{u(t)}_{\dot{W}^{1,1}(\bbS)}}\norm{u(t)}_{\dot{W}^{1,1}(\bbS)}.
$$
Such a nonlinear Sobolev inequality has the same flavour as the logarithmic Sobolev inequality \cite{gross1975logarithmic}
$$
\int_{\mathbb{R}}u^2\log(u^2)\,dx\leq \frac{1}{2}\log\left(\frac{2}{\pi e}\int_{\mathbb{R}}(\pa_xu)^2 \,dx\right).
$$
We will make use of this inequality in the study of the following one-dimensional \emph{arctan}-fast diffusion equation
\begin{subequations}\label{eq:1}
\begin{align}
\pt u-\px\arctan\left(\frac{\px u}{u}\right)&=0&& (x,t)\text{ on }\bbS\times[0,T],\\
u(x,0)&=u_0(x)&& x\text{ on }\bbS.
\end{align}
\end{subequations}
We will focus on the well-posedness for such equation and, in particular, we will establish the global existence of solutions for initial data satisfying certain properties.

There are several closely related problems in the literature. For instance, \eqref{eq:1} is a nonlinear diffusion somehow similar to the logarithmic fast diffusion equation
\begin{equation}\label{eq:log}
\pt u-\px\left(\frac{\px u}{u}\right)=0\qquad (x,t)\text{ on }\bbS\times[0,T].
\end{equation}
In fact, using
$$
\arctan(y)=y+\text{h.o.t.}
$$
equation \eqref{eq:log} can be obtained from \eqref{eq:1} when the effect of higher order nonlinearities is neglected.
Remarkably, \eqref{eq:log} is related to the Ricci flow and due to that this equation has been extensively studied in the past years by many authors (see \cite{vazquez1996fast} and the references therein). \\
We could also quote the fast diffusion equation
$$
\pt u=\px\left(\frac{\px u}{u^\sigma}\right)\q\t{with}\q\sigma>1.
$$
Equation \eqref{eq:1} is also related to the following nonlocal fast diffusion equation
\begin{equation}\label{eq:nonlocal}
\pt u-\px\arctan\left(\frac{-H u}{u}\right)=0\qquad (x,t)\text{ on }\bbS\times[0,T],
\end{equation}
where
$$
H u(x)=\frac{1}{2\pi}\pv \frac{u(y)}{\tan\left(\frac{x-y}{2}\right)}\,dy,
$$
is the Hilbert transform. Indeed, \eqref{eq:nonlocal} is obtained from \eqref{eq:1} by replacing the derivative inside the $\arctan(\cdot)$ by a term proportional to the Hilbert transform:
$$
\px\rightarrow -H.
$$
A similar procedure establishes relations between the KdV and the Benajmin-Ono equation or between the Sine-Gordon and the Sine-Hilbert equation \cite{degasperis}. Equation \eqref{eq:nonlocal} was derived by S. Steinerberger \cite{St2018} when studying how the distribution of roots behaves under iterated differentation (see also \cite{alazard2022dynamics,kiselev2022flow} for the mathematical study of some of its properties).\\
Actually, the principal reason that led us to study  problem \eqref{eq:1} is equation \eqref{eq:nonlocal}. Roughly speaking, our main aim is to understand how the nonlinear term driven by the $\arctan(\cdot)$ allows to prove the existence in the singular case, namely, where $u$ reaches zero.  Hence, a first step is to analyze what happens in the local case \eqref{eq:1}.

%yo diria algo como que la principal motivacion es la ecuacion no local que ponemos en la intro
%
%y que lo que queremos es entender bien como la nolinealidad atan permite demostrar existencia en el caso donde toca cero
%
%en el caso mas sencillo de un problema local

Finally, since previous equation \eqref{eq:1} can be written as
\begin{equation}\label{eq:patu}
\pat u -\frac{u\pa^2_x u-\pare{\pa_x u}^2}{u^2+\pare{\pa_xu}^2}=0,
\end{equation}
we observe that there is a striking similarity with the one dimensional relativistic heat equation
\[
\pat u=\pa_x\pare{\frac{u\pa_xu}{\sqrt{u^2+(\pa_xu)^2}}}=
\frac{u\pa^2_x u+\pare{\pa_xu}^2}{\sqrt{u^2+\pare{\pa_xu}^2}}
-u(\pa_xu)^2\frac{u+\pa_x^2u}{\pare{u^2+\pare{\pa_xu}^2}^{3/2}},
\]
see \cite{ACM,ACM2,ACM3,ACM4,ACMM,ACMM2,ACMS, ACMSV,CMSV}.

\subsection{Main results}

As anticipated, one of the main results of this paper is
\begin{theorem}\label{teo0}Let $0<u\in W^{1,1}(\bbS)$ be a function such that
$$
\|u\|_{L^1(\bbS)}=1.
$$
Then, the following nonlinear Sobolev inequalities hold true:
\begin{align*}
\int_{\bbS}\arctan\pare{\frac{\pa_x u}{u}}\pa_xu \,dx&\geq \arctan\pare{\norm{u}_{\dot{W}^{1,1}(\bbS)}}\norm{u}_{\dot{W}^{1,1}(\bbS)},\\
\int_{\bbS}\arctan\pare{\frac{\av{\pa_x u}}{u}}\frac{\av{\pa_xu}}{u}\,dx&\geq \frac{1}{4\pi}\arctan\pare{\frac{\|u\|_{\dot{W}^{1,1}(\bbS)}^2}{\||\pa_x u| u\|_{L^1(\bbS)}}}\frac{\|u\|_{\dot{W}^{1,1}(\bbS)}^2}{\||\pa_x u| u\|_{L^1(\bbS)}}.
\end{align*}
\end{theorem}
The assumption $\|u\|_{L^1(\bbS)}=1$ allows us to write the inequalities in Theorem \ref{teo0} in a simpler form. However, this hypothesis is not necessary and  an analogous result can be proved without this simplification.\\

We will make use of Theorem \ref{teo0} in the study \eqref{eq:1}, proving the following results.

\begin{theorem}\label{teo1}
Let $0< u_0\in H^3(\bbS)$ be the initial data. Then there exists a time $0<T\leq\infty$, $T=T(\|u_0\|_{H^3},\min_x u_0(x))$ and a unique positive solution to \eqref{eq:1}
$$
0<u\in C([0,T],H^3(\bbS)).
$$
\end{theorem}

\begin{theorem}\label{teo2}
Let $0< u_0\in H^3(\bbS)$ the initial data for \eqref{eq:1}. Then, as long as the unique positive solution $u$ to \eqref{eq:1} exists, the following properties hold:
\begin{itemize}
\item Maximum principle: $\|u(t)\|_{L^\infty(\bbS)}\leq\|u_0\|_{L^\infty(\bbS)},$
\item Mass conservation: $\|u(t)\|_{L^1(\bbS)}=\|u_0\|_{L^1(\bbS)},$
\item Entropy balance:
\begin{align*}%\label{eq:entropy}
\mathscr{H}(t)+\int_0^t\mathscr{D}(s) ds= \mathscr{H}(0),
\end{align*}
where
$$
\mathscr{H}(t)=\int_\bbS u(x,t)\log(u(x,t))-u(x,t)+1\,dx
$$
and
$$
\mathscr{D}(t)=\int_{\bbS}\arctan\pare{\frac{\pa_x u(x,t)}{u(x,t)}}\frac{\pa_xu(x,t)}{u(x,t)}\,dx.
$$
In particular, invoking Theorem \ref{teo0},
$$
\mathscr{H}(t)+\frac{1}{4\pi}\int_0^t \arctan\pare{\frac{\|u(s)\|_{\dot{W}^{1,1}(\bbS)}^2}{\||\pa_x u(s)| u(s)\|_{L^1(\bbS)}}}\frac{\|u(s)\|_{\dot{W}^{1,1}(\bbS)}^2}{\||\pa_x u(s)| u (s)\|_{L^1(\bbS)}}ds\leq \mathscr{H}(0).
$$
\item Energy balance:
\begin{equation*}%\label{identity:1}
\frac{1}{2}\|u(t)-\langle u_0\rangle\|_{L^2(\bbS)}^2+\int_0^t\mathcal{D}(s) ds= \frac{1}{2}\|u_0-\langle u_0\rangle\|_{L^2(\bbS)}^2,
\end{equation*}
where
$$
\mathcal{D}(t)=\int_{\bbS}\arctan\pare{\frac{\pa_x u}{u}}\pa_xu\,dx.
$$
In particular, invoking Theorem \ref{teo0},
\begin{equation*}%\label{identity:12}
\|u(t)-\langle u_0\rangle\|_{L^2(\bbS)}^2+\int_0^t\arctan\pare{\norm{u(s)}_{\dot{W}^{1,1}(\bbS)}}\norm{u(t)}_{\dot{W}^{1,1}(\bbS)} ds\leq  \|u_0-\langle u_0\rangle\|_{L^2(\bbS)}^2.
\end{equation*}
\item Energy decay:
$$
\|u(t)-\langle u_0\rangle\|_{L^2(\bbS)}\leq \|u_0-\langle u_0\rangle\|_{L^2(\bbS)}e^{-\frac{1}{2}\frac{\arctan (C\|u_0-\langle u_0\rangle\|_{L^2(\bbS)})}{\|u_0-\langle u_0\rangle\|_{L^2(\bbS)}}t}.
$$
\end{itemize}
\end{theorem}

We study further properties of solutions to \eqref{eq:1} passing to the following formulation of the problem.\\
We introduce the complex function $z(x,t)=u(x,t)+i\pa_xu(x,t)$ and, passing to polar coordinates, we rewrite \eqref{eq:1} as
\[
u(1+\tan^2 (\theta))\pat \theta=\pa_x^2\arctan\theta-\frac{\tan(\theta)}{1+\theta^2}\pa_x\theta=\frac{\pa_x^2\theta}{1+\theta^2}-\pare{  \frac{2\theta}{1+\theta^2}+\tan(\theta)}\frac{\pa_x\theta}{1+\theta^2},
\]
being
\[
\tan(\theta(x,t))=\frac{\pa_xu(x,t)}{u(x,t)}.
\]

\begin{theorem}\label{teo:prop-theta}
Let $0< u_0\in H^3(\bbS)$ the initial data for \eqref{eq:1}. Then, as long as the unique positive solution $u$ to \eqref{eq:1} exists, the following properties hold:
\begin{itemize}
\item Boundedness of the slope:
\[
\norm{\theta(t)}_{L^\infty(\bbS)}\le \norm{\theta(0)}_{L^\infty(\bbS)}.
\]
Furthermore,
\[
\norm{\pa_xu(t)}_{L^{\infty}(\bbS)} \le \norm{\pa_xu(0)}_{L^{\infty}(\bbS)}\frac{\max_x u_0(x)}{\min_xu_0(x)}.
\]
\item Lyapunov functional:
\begin{align*}
\frac{d}{dt}\int_{\bbS} u(1+\tan^2(\theta))\pare{\frac{\theta^2}{2}+\frac{\theta^4}{4}}\,dx&=-\int_{\bbS}
(\pa_x\theta)^2\,dx-\int_{\bbS}  \frac{(\pa_x\theta)^2(2+\theta^2)}{1+\theta^2}\theta\tan(\theta)
\,dx\\
&\q-\frac{1}{2}\int_{\bbS}  \frac{(\pa_x\theta)^2\theta^2(2+\theta^2)}{1+\theta^2}
(1+\tan^2(\theta))\,dx\\
&\le 0,
\end{align*}
i.e.
\begin{align*}
&\int_{\bbS} u(t)\pare{1+\pare{\frac{\pa_xu(t)}{u(t)}}^2}\pare{\arctan^2\pare{\frac{\pa_xu(t)}{u(t)}}+ \arctan^4\pare{\frac{\pa_xu(t)}{u(t)}}}\,dx\\
&\le 2\int_{\bbS} u_0\pare{1+\pare{\frac{\pa_xu(0)}{u_0}}^2}\pare{\arctan^2\pare{\frac{\pa_xu(0)}{u_0}}+\arctan^4\pare{\frac{\pa_xu(0)}{u_0}}}\,dx.
\end{align*}
\end{itemize}
In particular, we deduce that
 \[
u\in L^{\infty}(0,T;W^{1,\infty}(\bbS))\cap L^2(0,T;H^2(\bbS)).
\]
\end{theorem}

We also prove two different existence results to problem \eqref{eq:1} under lower regularity assumptions on the data.\\
In particular, we prove suitable uniform bounds to solutions of the approximating  problem
\begin{align*}
\pt u_\kappa-\px \arctan\left(\frac{\px u_\kappa}{u_\kappa}\right)&=0&&\t{in }\bbS\times(0,T_\kappa),\\
u_\kappa(x,0)&=\mathscr{J}_\kappa*u_0(x)&&\t{in }\bbS,
\end{align*}
 whose existence follows from  Theorem \ref{teo1}.\\

 We first exploit the Lyapunov functional found in Theorem \ref{teo:prop-theta} to obtain the following result.
\begin{theorem}[Existence through Lyapunov functional]\label{teo:lyap}
Let $u_0\in H^2(\bbS)$ with $u_0>0$, and $\norm{\arctan\pare{\frac{\pa_xu_0}{u_0}}}_{L^\infty(\bbS)}$ suitably small. Then, solutions to problem \eqref{eq:1} verify
\begin{align*}
&\int_{\bbS}\frac{\pare{u(t)\pa_x^2u(t)-(\pa_xu(t))^2}^2}{u(t)\pare{u^2(t)+(\pa_xu(t))^2}}\,dx\\
&\le \pare{\int_{\bbS}\frac{\pare{u_0\pa_x^2u(0)-(\pa_xu(0))^2}^2}{u_0\pare{u
_0^2+(\pa_xu(0))^2}}\,dx}e^{c\pare{\int_0^t\pare{\int_{\bbS}\frac{\pare{u(t)\pa_x^2u(t)-(\pa_xu(t))^2}^2}{\pare{u^2(t)+(\pa_xu(t))^2}^2}\,dx+1}\,dt}}.
\end{align*}
In particular,
 \[
u\in L^{\infty}(0,T;H^2(\bbS)).
\]
\end{theorem}

The next existence result we prove holds with Wiener data.

We recall that the definition of Wiener spaces is given by
\begin{align*}\label{Wienerhomo}
{{A}}^\alpha(\bbS)
&=\left\{u(x)\in L^1(\bbS):\q  \|u\|_{A^\alpha(\bbS)}:=\sum_{k\in\bbZ} |k|^\alpha|\widehat{u}(k)|<\infty\right\}.
\end{align*}
We will  make  use of the interpolation inequality \cite[Lemma 2.1]{BGB}:
\begin{equation*}%\label{interpol}
\|u\|_{A^p(\bbS)}\le \|u\|_{A^0(\bbS)(\bbS)}^{1-\theta}\|u\|_{A^q(\bbS)}^{\theta}\q\t{for}\q 0\le p\le q,\,\, \theta=\frac{p}{q}.
\end{equation*}

\begin{theorem}[Existence in Wiener spaces]\label{teo:wiener}
Let
\[
w(x,t)=\frac{u(x,t)-\langle u_0\rangle}{\langle u_0\rangle},
\]
with $u$ be the unique positive local solution to \eqref{eq:1}. Then, if $w_0=w(0)\in A^1(\bbS)$, and
\[
\norm{w_0}_{A^1(\bbS)}< \frac{1}{10},
\]
we have that
\[
w\in L^{\infty}(0,T;A^1(\bbS))\cap L^1(0,T; A^3(\bbS)).
\]
\end{theorem}
We observe that these two global existence results, although being stated in different spaces, have somehow the same flavour. Indeed, both have size restrictions in rather similar quantities, i.e.
$$
\frac{\pa_xu(0)}{u_0}\qquad\text{and}\qquad\frac{u(x,t)-\langle u_0\rangle}{\langle u_0\rangle}.
$$

\section*{Acknowledgments}
R.G-B thanks Javier G\'omez-Serrano for useful discussions. The work of M.M. was partially supported by Grant RYC2021-033698-I, funded by the Ministry of Science and Innovation/State Research Agency/10.13039/501100011033 and by the European Union "NextGenerationEU/Recovery, Transformation and Resilience Plan".\\
R.G-B and M.M. are funded by  the project "An\'alisis Matem\'atico Aplicado y Ecuaciones Diferenciales" Grant PID2022-141187NB-I00 funded by MCIN /AEI /10.13039/501100011033 / FEDER, UE and acronym "AMAED". This publication is part of the project PID2022-140494NA-I00 funded by MCIN/ AEI /10.13039/501100011033.

\section{Proof of Theorem \ref{teo0} - Sobolev inequalities}
First, we prove the following inequality
$$
\int_{\bbS}\arctan\pare{\frac{\av{\pa_x u}}{u}}\av{\pa_xu}\,dx\geq \arctan\pare{\norm{u}_{\dot{W}^{1,1}(\bbS)}}\norm{u}_{\dot{W}^{1,1}(\bbS)}.
$$
We observe that
$$
\int_{\bbS}\arctan\pare{\frac{\pa_x u}{u}}\pa_x u\,dx=\int_{\bbS}\arctan\pare{\frac{|\pa_x u|}{u}}|\pa_x u|\,dx
$$

We make use of the following inequality
\begin{equation*}
\arctan (z)\ge
\left\{
\begin{aligned}
&\frac{\arctan (\xi)}{\xi} z&&\qq\t{if }  z\le \xi,\\
&\arctan (\xi) &&\qq\t{if } z\ge \xi,
\end{aligned}
\right.
\end{equation*}
to deduce that
\begin{align*}
-\int_{\bbS}\arctan\pare{\frac{|\pa_x u|}{u}}|\pa_x u|\,dx&\le
%%
%&=-
-\left\{
\begin{aligned}
&\frac{\arctan \xi}{\xi} \intS \frac{|\pa_x u|^2}{u}\,dx&&\qq\t{if } \frac{|\pa_x u|}{u}\le \xi,\\
&\arctan \xi\intS \av{\pa_x u}\,dx &&\qq\t{if } \frac{|\pa_x u|}{u}\ge \xi.
\end{aligned}
\right.
\end{align*}
Since, by Hölder's inequality,
\[
\intS \av{\pa_x u}\,dx\le \pare{\intS \frac{(\pa_x u)^2}{u}\,dx}^\frac{1}{2}\pare{\intS u\,dx}^\frac{1}{2},
\]
we estimate
\begin{equation*}
 -\intS \frac{(\pa_x u)^2}{u}\,dx\le -\frac{\pare{\intS \av{\pa_x u}}^2}{\intS u\,dx},%\le -\frac{\intS \pare{\pa_x u}^2}{\intS u\,dx}
\end{equation*}
and we improve the bound as
\begin{equation}\label{eq:bound-atan}
-\int_{\bbS}\arctan\pare{\frac{|\pa_x u|}{u}}|\pa_x u|\,dx\le -
\left\{
\begin{aligned}
&\frac{\arctan \xi}{\xi} \frac{\pare{\intS \av{\pa_x u}}^2}{\intS u\,dx}&&\qq\t{if } \frac{|\pa_x u|}{u}\le \xi,\\
&\arctan \xi\intS \av{\pa_x u}\,dx &&\qq\t{if } \frac{|\pa_x u|}{u}\ge \xi.
\end{aligned}
\right.
\end{equation}

Let
\[
\xi= \intS \av{\pa_x u}\,dx
\]
in \eqref{eq:bound-atan}. This choice of $\xi$ leads to
\begin{equation*}%\label{eq:W11}
\begin{aligned}
&-\int_{\bbS}\arctan\pare{\frac{\av{\pa_x u}}{u}}\av{\pa_x u}\,dx\\
& \le-
\left\{
\begin{aligned}
 & \frac{\arctan\pare{\intS \av{\pa_x u}\,dx}}{\pare{\intS u\,dx}\intS \av{\pa_x u}\,dx}\pare{\intS \av{\pa_x u}}^2&\q\t{if } \frac{\av{\pa_x u}}{u}\le \intS \av{\pa_x u}\,dx,\\
&\arctan\pare{\intS \av{\pa_x u}\,dx}\pare{\intS \av{\pa_x u}\,dx}&\q\t{if } \frac{\av{\pa_x u}}{u}\ge\intS \av{\pa_x u}\,dx.
\end{aligned}
\right.
\end{aligned}
\end{equation*}
Then, using the simplification
$$
\int_\bbS u(x)\,dx=1,
$$
$$
-\int_{\bbS}\arctan\pare{\frac{\av{\pa_x u}}{u}}\av{\pa_x u}\,dx\le -\arctan\pare{\norm{u}_{\dot{W}^{1,1}(\bbS)}}\norm{u}_{\dot{W}^{1,1}(\bbS)}.
$$

Now, we prove the following inequality:
$$
\int_{\bbS}\arctan\pare{\frac{\av{\pa_x u}}{u}}\frac{\av{\pa_xu}}{u}\,dx\geq \frac{1}{4\pi}\arctan\pare{\frac{\|u\|_{\dot{W}^{1,1}(\bbS)}^2}{\||\pa_x u| u\|_{L^1(\bbS)}}}\frac{\|u\|_{\dot{W}^{1,1}(\bbS)}^2}{\||\pa_x u| u\|_{L^1(\bbS)}}.
$$
As before, we have that
\begin{align*}
\int_{\bbS}\arctan\pare{\frac{\pa_x u}{u}}\frac{\pa_x u}{u}\,dx&\geq
\left\{
\begin{aligned}
&\frac{\arctan \xi}{\xi} \intS \frac{(\pa_x u)^2}{u^2}\,dx&&\qq\t{if } \frac{\pa_x u}{u}\le \xi,\\
&\arctan \xi\intS \frac{\av{\pa_x u}}{u}\,dx &&\qq\t{if } \frac{\pa_x u}{u}\ge \xi.
\end{aligned}
\right.
\end{align*}
Furthermore,
\begin{equation*}
\intS \frac{|\pa_x u|}{u}\,dx\le \sqrt{2\pi}\pare{\intS \frac{|\pa_x u|^2}{u^2}\,dx}^{1/2},
\end{equation*}
We fix
$$
\xi=\intS \frac{\av{\pa_x u}}{u}\,dx,
$$
and we obtain
$$
\int_{\bbS}\arctan\pare{\frac{\av{\pa_x u}}{u}}\frac{\av{\pa_xu}}{u}\,dx\geq \frac{1}{4\pi}\arctan\pare{\intS \frac{\av{\pa_x u}}{u}\,dx}\intS \frac{\av{\pa_x u}}{u}\,dx.
$$
Using
\begin{equation*}
\intS \av{\pa_x u}\,dx\leq \pare{\intS \frac{\av{\pa_x u}}{u}\,dx}^{1/2}\pare{\intS \av{\pa_x u} u\,dx}^{1/2},
\end{equation*}
we conclude the desired inequality.

\section{Proof of Theorem \ref{teo1} - The existence result with $u_0\in H^3(\bbS)$}
\textbf{Well-posedness:} The well-posedness follows from the classical energy method \cite{Majda}, so we only sketch the proof.

We fix $\varepsilon, \kappa$ and $\delta$ three positive parameters and define the approximate problems
\begin{align*}
\pt u^{(\varepsilon,\kappa,\delta)}-\px\mathscr{J}_\kappa*\arctan\left(\frac{\px \mathscr{J}_\kappa*u^{(\varepsilon,\kappa,\delta)}}{\mathscr{J}_\kappa*u^{(\varepsilon,\kappa,\delta)}+\varepsilon}\right)&=\varepsilon \mathscr{J}_\kappa* \px^2 \mathscr{J}_\kappa*u^{(\varepsilon,\kappa,\delta)},\\
u^{(\varepsilon,\kappa,\delta)}(x,0)&=\mathscr{J}_\kappa*u^{(\varepsilon,\kappa,\delta)}_0(x)+\delta,
\end{align*}
where $\mathscr{J}_\kappa$ denotes the periodic heat kernel at time $\kappa$. The existence of a unique positive solution (up to a time $0<T(\varepsilon,\kappa,\delta)\leq\infty$) $u^{(\varepsilon,\kappa,\delta)}$ follows from an application of Picard Theorem in $H^3$ \cite{Majda}. Indeed, it is a tedious but straightforward computation to check that

$$
F^{(\varepsilon,\kappa,\delta)}(u^{(\varepsilon,\kappa,\delta)})=\varepsilon \mathscr{J}_\kappa* \px^2 \mathscr{J}_\kappa*u^{(\varepsilon,\kappa,\delta)}+\px\mathscr{J}_\kappa*\arctan\left(\frac{\px \mathscr{J}_\kappa*u^{(\varepsilon,\kappa,\delta)}}{\mathscr{J}_\kappa*u^{(\varepsilon,\kappa,\delta)}+\varepsilon}\right)
$$
satisfies
$$
F^{(\varepsilon,\kappa,\delta)}:H^3\rightarrow H^3
$$
and it is a Lipschitz operator
$$
\|F^{(\varepsilon,\kappa,\delta)}(u^{(\varepsilon,\kappa,\delta)})-F^{(\varepsilon,\kappa,\delta)}(v^{(\varepsilon,\kappa,\delta)})\|_{H^3}\leq C(\varepsilon,\kappa,\delta)\|u^{(\varepsilon,\kappa,\delta)}-v^{(\varepsilon,\kappa,\delta)}\|_{H^3}.
$$
As a consequence, it exists a positive time of existence $0<T(\varepsilon,\kappa,\delta)$ and a smooth solution
$$
u^{(\varepsilon,\kappa,\delta)}\in C([0,T(\varepsilon,\kappa,\delta)],H^3(\bbS)).
$$
The next step is to obtain uniform estimates for $0<T^*<T(\varepsilon,\kappa,\delta)$.

It is easy to obtain $\kappa-$uniform bounds. Indeed, we test against
$$
u^{(\varepsilon,\kappa,\delta)}
$$
and, using the properties of the convolution, the symmetry of the heat kernel and Young's inequality, obtain
$$
\frac{d}{dt}\|u^{(\varepsilon,\kappa,\delta)}(t)\|_{L^2(\bbS)}^2+\varepsilon\|\mathscr{J}_\kappa*u^{(\varepsilon,\kappa,\delta)}(t)\|_{\dot{H}^1(\bbS)}^2\leq C(\varepsilon).
$$
The higher order estimate can be obtained in a similar fashion. We test against
$$
-\px^6 u^{(\varepsilon,\kappa,\delta)}
$$
and integrate by parts. We find
$$
\frac{1}{2}\frac{d}{dt}\|u^{(\varepsilon,\kappa,\delta)}(t)\|_{\dot{H}^3(\bbS)}^2+\varepsilon\|\mathscr{J}_\kappa*u^{(\varepsilon,\kappa,\delta)}(t)\|_{\dot{H}^4(\bbS)}^2=-\int_{\bbS}\px^3\arctan\left(\frac{\px \mathscr{J}_\kappa*u^{(\varepsilon,\kappa,\delta)}}{\mathscr{J}_\kappa*u^{(\varepsilon,\kappa,\delta)}+\varepsilon}\right)\px^4\mathscr{J}_\kappa*u^{(\varepsilon,\kappa,\delta)}\,dx.
$$
The nonlinear term can be handled easily using the parabolicity and integration by parts. Then we conclude
$$
\frac{1}{2}\frac{d}{dt}\|u^{(\varepsilon,\kappa,\delta)}(t)\|_{\dot{H}^3(\bbS)}^2\leq \mathcal{P}(\|\mathscr{J}_\kappa*u^{(\varepsilon,\kappa,\delta)}(t)\|_{\dot{H}^3(\bbS)}^2)\leq \mathcal{P}(\|u^{(\varepsilon,\kappa,\delta)}(t)\|_{\dot{H}^3(\bbS)}^2),
$$
which leads to the desired $\kappa-$uniform estimate. We can pass to the limit and obtain a solution
$$
u^{(\varepsilon,\delta)}\in C([0,T(\varepsilon,\delta)],H^3(\bbS)).
$$
In order to obtain $\varepsilon-$uniform estimates we define
\begin{equation*}%\label{ineq:0}
\mathscr{E}^{(\varepsilon,\delta)}(t)=\frac{1}{\min_x u^{(\varepsilon,\delta)}(x,t)}+\|u^{(\varepsilon,\delta)}(t)\|_{H^3(\bbS)}.
\end{equation*}
Now we use a pointwise argument (see  \cite{cordobacordoba} for more details). Being continuous and the domain a compact set, the solution has at least a minimum:
$$
m^{(\varepsilon,\delta)}(t)=\min_ x u^{(\varepsilon,\delta)}(x,t)=u^{(\varepsilon,\delta)}(\underline{x}_t,t).
$$
Because of the positivity of the initial data, we have that $m^{(\varepsilon,\delta)}(0)>\delta>0$. Following the argument in \cite{cordobacordoba}, we have that
$$
\frac{d}{dt}m^{(\varepsilon,\delta)}(t)=\pt u^{(\varepsilon,\delta)}(\underline{x}_t,t)=\frac{\px^2 u^{(\varepsilon,\delta)}(\underline{x}_t)}{m^{(\varepsilon,\delta)}(t)}\;\text{ a.e.}.
$$
Indeed, due to the smoothness of $u^{(\varepsilon,\delta)}$ in space and time we have that $m^{(\varepsilon,\delta)}(t)$ is Lipschitz. To see that we use the reverse triangle inequality and find that
\begin{equation*}
|\min_x u^{(\varepsilon,\delta)}(t_1,x)-\min_x u^{(\varepsilon,\delta)}(t_2,x)|\leq \min_x (|u^{(\varepsilon,\delta)}(t_1,x)-u^{(\varepsilon,\delta)}(t_2,x)|)\leq C|t_1-t_2|.
\end{equation*}
Using Rademacher's Theorem, we have that $\min_x u^{(\varepsilon,\delta)}(t,x)$ is differentiable almost everywhere. Thus, using that $x_t$ is the point of minimum, we get that
\begin{eqnarray*}\frac{d}{dt}m^{(\varepsilon,\delta)}(t) &=& \lim_{h_j\rightarrow 0} \frac{m^{(\varepsilon,\delta)}(t+h_j)-m^{(\varepsilon,\delta)}(t)}{h_j} \\
&=& \lim_{h_j\rightarrow 0}\frac{u^{(\varepsilon,\delta)}(x_{t+h_j},t+h_j) - u^{(\varepsilon,\delta)}(x_t,t)}{h_j}\\
&\geq& \pat u^{(\varepsilon,\delta)}(x_t,t).
\end{eqnarray*}
In the same way we compute
\begin{eqnarray*}\frac{d}{dt}m^{(\varepsilon,\delta)}(t) &\leq& \pat u^{(\varepsilon,\delta)}(x_t,t).
\end{eqnarray*}

Then,
\begin{equation*}%\label{ineq:1}
\frac{d}{dt}\frac{1}{\min_x u(x,t)}=-\frac{\pt u(\underline{x}_t,t)}{m(t)^2}\leq C\frac{\|u(t)\|_{H^3(\bbS)}}{m(t)^3}\leq C(\mathscr{E}(t))^4.
\end{equation*}
We can estimate the evolution of the $H^3$ norm with the previous ideas together with the fact that
$$
\frac{1}{u^{(\varepsilon,\delta)}+\varepsilon}\leq\frac{1}{u^{(\varepsilon,\delta)}}\leq \frac{1}{m^{(\varepsilon,\delta)}(t)}.
$$
Thus, finally we conclude
$$
\frac{d}{dt}\mathscr{E}^{(\varepsilon,\delta)}(t)\leq \mathcal{P}(\mathscr{E}^{(\varepsilon,\delta)}(t)).
$$
Then, we obtain a $\varepsilon-$uniform bound and we can pass to the limit. We find
$$
u^{(\delta)}\in C([0,T(\delta)],H^3(\bbS)).
$$
With the previous ideas we can also pass to the limit in $\delta$ and we conclude the local existence of classical solution
$$
u\in C([0,T],H^3(\bbS)).
$$
To obtain the uniqueness we proceed with a standard contradiction argument.

\section{Proof of Theorem \ref{teo2} - Properties of the solution}
\textbf{Maximum principle:} Equation \eqref{eq:1} can be written as
\[
\pat u -\frac{u\pa^2_x u-\pare{\pa_x u}^2}{u^2+\pare{\pa_xu}^2}=0.
\]
Using the pointwise method \cite{cordobacordoba}, we find that
$$
M(t)=\max_x u(x,t)=u(\overline{x}_t,t),
$$
and
$$
m(t)=\min_x u(x,t)=u(\underline{x}_t,t),
$$
satisfy
$$
\frac{d}{dt}M\leq 0, \text{ a.e.}\qquad\t{and}\qq \frac{d}{dt}m\geq 0, \text{ a.e.}.
$$
Integrating in time we conclude this part. We observe that, in particular
$$
\|u(t)\|_{L^\infty(\bbS)}\leq\|u_0\|_{L^\infty(\bbS)}.
$$

\medskip

\textbf{Conservation of mass:} The conservation of mass follows from the sign propagation and an integration in space.

\medskip

\textbf{Entropy balance:} The evolution of the entropy can be easily computed and we find that
$$
\frac{d}{dt}\mathscr{H}(t)+\int_\bbS \arctan\pare{\frac{\partial_x u}{u}}\frac{\partial_x u}{u} \,dx=0.
$$
Using that
$$
\arctan(z)z=\arctan(|z|)|z|\geq0,
$$
we conclude. Using Theorem \ref{teo0} we conclude the desired estimate.

\medskip

\textbf{Energy balance:} The evolution of the $L^2$ energy can be computed similarly. We observe that the mean is preserved (see above). To estimate the decay, we compute the following:
$$
\frac{1}{2}\frac{d}{dt}\|u(t)-\langle u(t)\rangle\|_{L^2(\bbS)}^2=\int_\bbS \partial_t u(u-\langle u_0\rangle)\,dx=\int_\bbS \partial_t u u\,dx.
$$
Integrating by parts we find
$$
\mathcal{D}(t)=\int_{\bbS}\arctan\pare{\frac{\pa_x u}{u}}\pa_xu\,dx\geq0.
$$
In particular, invoking Theorem \ref{teo0},
\begin{equation*}%\label{identity:1a}
\|u(t)-\langle u_0\rangle\|_{L^2(\bbS)}^2+\int_0^t\arctan\pare{\norm{u(s)}_{\dot{W}^{1,1}(\bbS)}}\norm{u(s)}_{\dot{W}^{1,1}(\bbS)} ds\leq  \|u_0-\langle u_0\rangle\|_{L^2(\bbS)}^2,
\end{equation*}
and we obtain that
$$
\|u(t)-\langle u_0\rangle\|_{L^2(\bbS)}
$$
decays.

\medskip

\textbf{Energy decay:} Using that, for certain point $\tilde{x}_t$,
$$
u(x,t)-\langle u_0\rangle=u(x,t)-u(\tilde{x}_t,t)=\int_{\tilde{x}_t}^x\partial_x u(y,t)\,dy,
$$
we conclude the Poincar\'e inequality
$$
C\|u(t)-\langle u_0\rangle\|_{L^2(\bbS)}\leq \|\partial_x u(t)\|_{L^1(\bbS)}.
$$
As a consequence,
$$
\arctan\pare{\norm{u(t)}_{\dot{W}^{1,1}(\bbS)}}\norm{u(t)}_{\dot{W}^{1,1}(\bbS)}\geq\arctan\pare{C\|u(t)-\langle u_0\rangle\|_{L^2(\bbS)}}C\|u(t)-\langle u_0\rangle\|_{L^2(\bbS)}.
$$
Thus, integrating by parts and estimating using Theorem \ref{teo0},
\begin{equation*}%\label{identity:1b}
2\frac{d}{dt}\|u(t)-\langle u_0\rangle\|_{L^2(\bbS)}+C\arctan\pare{C\|u(t)-\langle u_0\rangle\|_{L^2(\bbS)}}\leq 0.
\end{equation*}
Using that
\begin{equation*}
\arctan (z)\ge
\frac{\arctan (C\|u_0-\langle u_0\rangle\|_{L^2(\bbS)})}{C\|u_0-\langle u_0\rangle\|_{L^2(\bbS)}} z\qq \text{if}\qq  z\le C\|u_0-\langle u_0\rangle\|_{L^2(\bbS)},
\end{equation*}
the previous decay of the $L^2$ energy translates into the following inequality
\begin{equation*}%\label{identity:1c}
2\frac{d}{dt}\|u(t)-\langle u_0\rangle\|_{L^2(\bbS)}+\frac{\arctan \pare{C\|u_0-\langle u_0\rangle\|_{L^2(\bbS)}}}{\|u_0-\langle u_0\rangle\|_{L^2(\bbS)}}\|u(t)-\langle u_0\rangle\|_{L^2(\bbS)}\leq 0,
\end{equation*}
from where we conclude the
$$
\|u(t)-\langle u_0\rangle\|_{L^2(\bbS)}\leq \|u_0-\langle u_0\rangle\|_{L^2(\bbS)}e^{-\frac{1}{2}\frac{\arctan \pare{C\|u_0-\langle u_0\rangle\|_{L^2(\bbS)}}}{\|u_0-\langle u_0\rangle\|_{L^2(\bbS)}}t}.
$$

\section{Proof of Theorem \ref{teo:prop-theta} - Further properties of the solution with the $\theta$ formulation}
We consider the complex value $z(x,t)=u(x,t)+i\pa_xu(x,t)$. If we write this quantity in polar coordinates, we know that
\begin{equation}\label{eq:arg}
\tan(\theta(x,t))=\frac{\pa_xu(x,t)}{u(x,t)},\q\t{hence}\q \theta(x,t)=\arctan\pare{\frac{\pa_xu(x,t)}{u(x,t)}}.
\end{equation}
Then, we use \eqref{eq:1} to deduce the evolution equation of $\theta$:
\begin{equation}\label{eq:theta}
u(1+\tan^2 (\theta))\pat \theta=\pa_x^2\arctan\theta-\frac{\tan(\theta)}{1+\theta^2}\pa_x\theta=\frac{\pa_x^2\theta}{1+\theta^2}-\pare{  \frac{2\theta}{1+\theta^2}+\tan(\theta)}\frac{\pa_x\theta}{1+\theta^2}.
\end{equation}

\textbf{Boundedness of the slope:} Taking in mind \eqref{eq:arg} and \eqref{eq:theta}, let us take
\[
\Psi(t)=\max_x\theta(x,t)=\theta(x_\Psi,t).
\]
Then, a.e. $t\in[0,T]$, we have that
\[
u(x_\Psi,t)(1+\tan^2 \Psi(t))\pat \Psi(t)=\frac{\pa_x^2\theta(x_\Psi,t)}{1+\Psi^2(t)}\le0.
\]
Observing also that $u(x_\Psi,t)(1+\tan^2 \Psi(t))>0$, we obtain that
\[
\max_x\theta(x,t)=\Psi(t)\le \Psi(0)=\max_x\theta(x,0).
\]
We can repeat the same argument for
\[
\Phi(t)=\min_x\theta(x,t)=\theta(x_\Phi,t),
\]
getting that
\[
\min_x\theta(x,t)=\Phi(t)\ge \Phi(0)=\min_x\theta(x,0).
\]
This implies that the function $\theta$ is bounded in the $x$ variable a.e. $t$, i.e.
\[
\norm{\theta(t)}_{L^\infty(\bbS)}\le \norm{\theta(0)}_{L^\infty(\bbS)}\qq \t{a.e. }t\in[0,T].
\]
We now recall the definition of $\theta$ to deduce that
\begin{align*}
\norm{\frac{\pa_xu(t)}{u(t)} }_{L^\infty(\bbS)}
&= \norm{\tan(\theta(t)) }_{L^\infty(\bbS)} \le  \norm{\tan(\theta(0)) }_{L^\infty(\bbS)} ,
\end{align*}
and the desired estimate follows from
  the boundedness of $u$:
\begin{align*}
\norm{\pa_xu(t)}_{L^{\infty}(\bbS)}&=\norm{\frac{\pa_xu(t)}{u(t)}u(t) }_{L^\infty(\bbS)}\le  \norm{\frac{\pa_xu(t)}{u(t)} }_{L^\infty(\bbS)}\norm{u(t) }_{L^\infty(\bbS)}\le  \norm{\tan(\theta(0)) }_{L^\infty(\bbS)}\norm{u_0 }_{L^\infty(\bbS)}\\
&\le  \norm{\pa_xu(0)}_{L^{\infty}(\bbS)}\frac{\max_x u_0(x)}{\min_xu_0(x)}.
\end{align*}

\medskip

\textbf{Lyapunov functional:}
We are going to prove that the functional
\[
L(u)=\int_{\bbS} u(1+\tan^2(\theta))\pare{\frac{\theta^2}{2}+\frac{\theta^4}{4}}\,dx
\]
is a Lyapunov functional. To this aim, we compute
\begin{align}
\frac{d}{dt}L(u)&=\int_{\bbS} \pat u(1+\tan^2(\theta))\pare{\frac{\theta^2}{2}+\frac{\theta^4}{4}}\,dx +\int_{\bbS} u\pat(1+\tan^2(\theta))\pare{\frac{\theta^2}{2}+\frac{\theta^4}{4}}\,dx\nonumber\\
&\q+\int_{\bbS} u(1+\tan^2(\theta))\pat\pare{\frac{\theta^2}{2}+\frac{\theta^4}{4}}\,dx\nonumber\\
&=\int_{\bbS} \pat u(1+\tan^2(\theta))\pare{\frac{\theta^2}{2}+\frac{\theta^4}{4}}\,dx +\int_{\bbS} u\tan(\theta)(1+\tan^2(\theta))\pare{\theta^2+\frac{\theta^4}{2}}\pat\theta\,dx\nonumber\\
&\q+\int_{\bbS} u(1+\tan^2(\theta))\pare{\theta+\theta^3}\pat\theta\,dx,\label{eq:lyap1}
\end{align}
and we use the equations
\begin{align*}
u(1+\tan^2 (\theta))\pat \theta&=\frac{\pa_x^2\theta}{1+\theta^2}-\pare{  \frac{2\theta}{1+\theta^2}+\tan(\theta)}\frac{\pa_x\theta}{1+\theta^2},\\
\pat u&=\frac{u\pa^2_x u-\pare{\pa_x u}^2}{u^2+\pare{\pa_xu}^2}=\pa_x\theta,
\end{align*}
to rewrite \eqref{eq:lyap1} as
\begin{align}
&\frac{d}{dt}L(u)\nonumber\\
&=\int_{\bbS} \pa_x\theta(1+\tan^2(\theta))\pare{\frac{\theta^2}{2}+\frac{\theta^4}{4}}\,dx
+\int_{\bbS}\theta u(1+\tan^2(\theta)) \bra{1+\theta^2+ \tan(\theta)\pare{\theta+\frac{\theta^3}{2}}}\pat\theta\,dx\nonumber\\
&=\int_{\bbS} \pa_x\theta(1+\tan^2(\theta))\pare{\frac{\theta^2}{2}+\frac{\theta^4}{4}}\,dx
\nonumber\\
&\q+\int_{\bbS}\theta  \bra{1+\theta^2+ \tan(\theta)\pare{\theta+\frac{\theta^3}{2}}}\pare{\frac{\pa_x^2\theta}{1+\theta^2}-\pare{  \frac{2\theta}{1+\theta^2}+\tan(\theta)}\frac{\pa_x\theta}{1+\theta^2}}\,dx.\label{eq:lyap2}
\end{align}
We claim that \eqref{eq:lyap2} is equivalent to
\begin{align}
\frac{d}{dt}\int_{\bbS} u(1+\tan^2(\theta))\pare{\frac{\theta^2}{2}+\frac{\theta^4}{4}}\,dx=
\int_{\bbS}\theta  \bra{1+\theta^2+ \tan(\theta)\pare{\theta+\frac{\theta^3}{2}}}\frac{\pa_x^2\theta}{1+\theta^2}\,dx.\label{eq:lyap3}
\end{align}
Indeed, the integrals
\[
\int_{\bbS} \pa_x\theta(1+\tan^2(\theta))\pare{\frac{\theta^2}{2}+\frac{\theta^4}{4}}\,dx
\]
and
\[
\int_{\bbS}\theta  \bra{1+\theta^2+ \tan(\theta)\pare{\theta+\frac{\theta^3}{2}}}\pare{ \frac{2\theta}{1+\theta^2}+\tan(\theta)}\frac{\pa_x\theta}{1+\theta^2}\,dx,
\]
can be rewritten both in the form
\[
\int_{\bbS}\phi'(\theta)\pa_x\theta\,dx,
\]
being
\[
\phi(\theta)=\int_{0}^{\theta}\psi(y)\,dy,
\]
with $\psi(y)$ either
\[
\psi(y)=(1+\tan^2(y))\pare{\frac{y^2}{2}+\frac{y^4}{4}},
\]
or
\[
\psi(y)=  \bra{1+y^2+ \tan (y)\pare{y+\frac{y^3}{2}}}\pare{ \frac{2y}{1+y^2}+\tan (y)}\frac{y}{1+y^2}.
\]
Then, by periodic boundary conditions, we obtain that
\[
\int_{\bbS}\phi'(\theta)\pa_x\theta\,dx=\phi(\theta)\biggl|_{-\pi}^\pi=0.
\]
We come back to \eqref{eq:lyap3}, and we rewrite the r.h.s. as
\begin{align*}
\frac{d}{dt}\int_{\bbS} u(1+\tan^2(\theta))\pare{\frac{\theta^2}{2}+\frac{\theta^4}{4}}\,dx&=
\int_{\bbS}
\theta\pa_x^2\theta\,dx
+
\frac{1}{2}\int_{\bbS}
\theta^2\tan(\theta)\pa_x^2\theta\,dx\\
&\q
+\frac{1}{2}\int_{\bbS} \tan(\theta)\pa_x^2\theta\,dx-\frac{1}{2}\int_{\bbS} \frac{\tan(\theta)}{1+\theta^2}\pa_x^2\theta\,dx.
\end{align*}
We compute the above integrals by integration by parts:
\begin{align*}
\int_{\bbS}
\theta\pa_x^2\theta\,dx&=-\int_{\bbS}
(\pa_x\theta)^2\,dx,
\\
\frac{1}{2}\int_{\bbS}
\theta^2\tan(\theta)\pa_x^2\theta\,dx&=-\frac{1}{2}\int_{\bbS}  (\pa_x\theta)^2
\pare{2\theta\tan(\theta)+\theta^2(1+\tan^2(\theta))}\,dx,
\\
\frac{1}{2}\int_{\bbS} \tan(\theta)\pa_x^2\theta\,dx&=-\frac{1}{2}\int_{\bbS}  (\pa_x\theta)^2
(1+\tan^2(\theta))\,dx,
\\
-\frac{1}{2}\int_{\bbS} \frac{\tan(\theta)}{1+\theta^2}\pa_x^2\theta\,dx&=
\frac{1}{2}\int_{\bbS}  \frac{(\pa_x\theta)^2}{1+\theta^2}
(1+\tan^2(\theta))\,dx-
\int_{\bbS}  \frac{(\pa_x\theta)^2\theta\tan(\theta)}{1+\theta^2}
\,dx.
\end{align*}
Summing up all the above integrals, we find that the r.h.s. of \eqref{eq:lyap3} simplifies as
\begin{align*}
\frac{d}{dt}\int_{\bbS} u(1+\tan^2(\theta))\pare{\frac{\theta^2}{2}+\frac{\theta^4}{4}}\,dx&=-\int_{\bbS}
(\pa_x\theta)^2\,dx-\int_{\bbS}  \frac{(\pa_x\theta)^2(2+\theta^2)}{1+\theta^2}\theta\tan(\theta)
\,dx\\
&\q-\frac{1}{2}\int_{\bbS}  \frac{(\pa_x\theta)^2\theta^2(2+\theta^2)}{1+\theta^2}
(1+\tan^2(\theta))\,dx.
\end{align*}
Using the properties of $\tan$, we have that $x\tan(x)\ge0$, and hence we conclude that
\[
\frac{d}{dt}\int_{\bbS} u(1+\tan^2(\theta))\pare{\frac{\theta^2}{2}+\frac{\theta^4}{4}}\,dx\le 0,
\]
and
$$
\int_0^t\int_{\bbS}(\partial_x \theta)^2\,dx\leq L(u_0).
$$

\section{Proof of Theorem \ref{teo:lyap} - Existence through Lyapunov functional}
Let $-\pa_{x}^2\theta$ be the test function in \eqref{eq:theta}. Then, the integral in space reads
\begin{align}
-\int_{\bbS}u\pare{1+\tan^2 \pare{\theta}}\pat \theta\pa_{x}^2\theta\,dx&=
-\int_{\bbS}\frac{\pare{\pa_x^2\theta}^2}{1+\theta^2}\,dx
+2\int_{\bbS}\frac{\theta\pa_x\theta\pa_{x}^2\theta}{\pare{1+\theta^2}^2}\,dx\nonumber\\
&\q
+\int_{\bbS}\tan\pare{\theta}\frac{\pa_x\theta\pa_{x}^2\theta}{1+\theta^2}\,dx.\label{eq:thetaH1}
\end{align}
We integrate by parts the integral in the l.h.s.:
\begin{align*}
&-\int_{\bbS}u\pare{1+\tan^2 \pare{\theta}}\pat \theta\pa_{x}^2\theta\,dx\\
&=\int_{\bbS}u\pare{1+\tan^2 \pare{\theta}}\pa_{x}\theta\pa_{tx} \theta\,dx
+\int_{\bbS}\pa_x\pare{u\pare{1+\tan^2 \pare{\theta}}}\pa_{x}\theta\pa_{t} \theta\,dx\\
&=\frac{1}{2}\frac{d}{dt}\int_{\bbS}u\pare{1+\tan^2 \pare{\theta}}\pare{\pa_{x}\theta}^2\,dx
-I_1+I_2,
\end{align*}
with
\begin{align*}
I_1&=\frac{1}{2}\int_{\bbS}\pat\pare{u\pare{1+\tan^2 \pare{\theta}}}\pare{\pa_{x}\theta}^2 \,dx,\\
I_2&=\int_{\bbS}\pa_x\pare{u\pare{1+\tan^2 \pare{\theta}}}\pa_{x}\theta\pa_{t} \theta\,dx.
\end{align*}
We use the equations
\begin{align*}
u(1+\tan^2 (\theta))\pat \theta&=\frac{\pa_x^2\theta}{1+\theta^2}-\pare{  \frac{2\theta}{1+\theta^2}+\tan(\theta)}\frac{\pa_x\theta}{1+\theta^2},\\
\pat u&=\frac{u\pa^2_x u-\pare{\pa_x u}^2}{u^2+\pare{\pa_xu}^2}=\pa_x\theta,\\
\frac{\pa_x u}{u}&=\tan(\theta),
\end{align*}
to rewrite $I_1$ and $I_2$ as follows:
\begin{align*}
I_1&=\frac{1}{2}\int_{\bbS}\pat u \pare{1+\tan^2 \pare{\theta}}\pare{\pa_{x}\theta}^2 \,dx+ \int_{\bbS} \tan(\theta)u\pare{1+\tan^2 \pare{\theta}}\pat \theta\pare{\pa_{x}\theta}^2 \,dx\\
&=\frac{1}{2}\int_{\bbS}   \pare{1+  \tan^2 \pare{\theta}}\pare{\pa_{x}\theta}^3 \,dx
+\int_{\bbS}\tan(\theta)\pare{\pa_{x}\theta}^2\frac{\pa_x^2\theta}{1+\theta^2}\,dx\\
&\q -2\int_{\bbS}\tan(\theta) \frac{\theta\pare{\pa_{x}\theta}^3}{(1+\theta^2)^2}\,dx-\int_{\bbS}\tan^2(\theta)\frac{\pare{\pa_{x}\theta}^3}{1+\theta^2}\,dx,
\\
I_2&=\int_{\bbS}\tan(\theta) u\pare{1+\tan^2 \pare{\theta}}\pa_{x}\theta(1+2\pa_x\theta)\pa_{t} \theta\,dx\\
&=\int_{\bbS}\tan(\theta)  \frac{\pa_x\theta\pa_x^2\theta}{1+\theta^2} \,dx-2\int_{\bbS}\tan(\theta)  \frac{(\pa_x\theta)^2 \theta}{(1+\theta^2)^2} \,dx-\int_{\bbS}\tan^2(\theta)  \frac{(\pa_x\theta)^2}{1+\theta^2} \,dx\\
&\q+2\int_{\bbS}\tan(\theta)  \frac{(\pa_x\theta)^2\pa_x^2\theta}{1+\theta^2} \,dx-4\int_{\bbS}\tan(\theta)  \frac{(\pa_x\theta)^3 \theta}{(1+\theta^2)^2} \,dx-2\int_{\bbS}\tan^2(\theta)  \frac{(\pa_x\theta)^3}{1+\theta^2} \,dx.
\end{align*}
The difference among $I_1$ and $I_2$ reads
\begin{align*}
I_1-I_2&=-\int_{\bbS}\tan(\theta)  \frac{\pa_x\theta\pa_x^2\theta}{1+\theta^2} \,dx+2\int_{\bbS}\tan(\theta)  \frac{(\pa_x\theta)^2 \theta}{(1+\theta^2)^2} \,dx+\int_{\bbS}\tan^2(\theta)  \frac{(\pa_x\theta)^2}{1+\theta^2} \,dx\\
&\q-\int_{\bbS}\tan(\theta)  \frac{(\pa_x\theta)^2\pa_x^2\theta}{1+\theta^2} \,dx+2\int_{\bbS}\tan(\theta)  \frac{(\pa_x\theta)^3 \theta}{(1+\theta^2)^2} \,dx+\int_{\bbS}\tan^2(\theta)  \frac{(\pa_x\theta)^3}{1+\theta^2} \,dx\\
&\q+\frac{1}{2}\int_{\bbS}   \pare{1+  \tan^2 \pare{\theta}}\pare{\pa_{x}\theta}^3 \,dx .
\end{align*}
Then, \eqref{eq:thetaH1} becomes
\begin{align}
&\frac{1}{2}\frac{d}{dt}\int_{\bbS}u\pare{1+\tan^2 \pare{\theta}}\pare{\pa_{x}\theta}^2\,dx+\int_{\bbS}\frac{\pare{\pa_x^2\theta}^2}{1+\theta^2}\,dx\nonumber\\
&=I_1-I_2
+2\int_{\bbS}\frac{\theta\pa_x\theta\pa_{x}^2\theta}{\pare{1+\theta^2}^2}\,dx
+\int_{\bbS}\tan\pare{\theta}\frac{\pa_x\theta\pa_{x}^2\theta}{1+\theta^2}\,dx\nonumber\\
&= 2\int_{\bbS}\tan(\theta)  \frac{(\pa_x\theta)^2 \theta}{(1+\theta^2)^2} \,dx+\int_{\bbS}\tan^2(\theta)  \frac{(\pa_x\theta)^2}{1+\theta^2} \,dx-\int_{\bbS}\tan(\theta)  \frac{(\pa_x\theta)^2\pa_x^2\theta}{1+\theta^2} \,dx\nonumber\\
&\q+\int_{\bbS}\tan^2(\theta)  \frac{(\pa_x\theta)^3}{1+\theta^2} \,dx +\frac{1}{2}\int_{\bbS}(1+\tan^2(\theta) )  \pare{\pa_{x}\theta}^3 \,dx
+2\int_{\bbS}\frac{\theta\pa_x\theta\pa_{x}^2\theta}{\pare{1+\theta^2}^2}\,dx . \label{eq:thetaH1-2}
\end{align}
We easily estimate the sum of the first two terms in \eqref{eq:thetaH1-2}
\[
2\int_{\bbS}\tan(\theta)  \frac{(\pa_x\theta)^2 \theta}{(1+\theta^2)^2} \,dx+\int_{\bbS}\tan^2(\theta)  \frac{(\pa_x\theta)^2}{1+\theta^2} \,dx
\]
using that
\[
\pa_x\theta\in L^2(0,T;L^2(\bbS)),
\]
and
\[
2\frac{\tan(\theta(t,x))\theta(t,x)}{1+\theta^2(t,x)}+\tan^2(\theta(t,x))\le \norm{2\tan(\theta(0))\theta(0)+\tan^2(\theta(0))}_{L^\infty(\bbS)}\le c.
\]
Furthermore, we can absorb the last term in \eqref{eq:thetaH1-2} using Young's inequality and the boundedness in time and space of $\theta$:
\[
2\int_{\bbS}\frac{\theta\pa_x\theta\pa_{x}^2\theta}{\pare{1+\theta^2}^2}\,dx\le \eps\int_{\bbS}\frac{\pare{\pa_x^2\theta}^2}{1+\theta^2}\,dx+c(\eps)\int_{\bbS}\frac{\theta^2(\pa_x\theta)^2}{\pare{1+\theta^2}^3}\,dx\le \eps\int_{\bbS}\frac{\pare{\pa_x^2\theta}^2}{1+\theta^2}\,dx+c(\eps)\int_{\bbS}(\pa_x\theta)^2\,dx.
\]
So far, \eqref{eq:thetaH1-2} can be bounded with
\begin{align}
&\frac{1}{2}\frac{d}{dt}\int_{\bbS}u\pare{1+\tan^2 \pare{\theta}}\pare{\pa_{x}\theta}^2\,dx+(1-\eps)\int_{\bbS}\frac{\pare{\pa_x^2\theta}^2}{1+\theta^2}\,dx\nonumber\\
&\le c(\eps)\int_{\bbS}(\pa_x\theta)^2\,dx+\int_{\bbS}\av{\tan(\theta) } \frac{(\pa_x\theta)^2\av{\pa_x^2\theta}}{1+\theta^2} \,dx +\frac{1}{2}\int_{\bbS}\pare{
1+\frac{3+\theta^2}{1+\theta^2}\tan^2(\theta)
} \pare{\pa_{x}\theta}^3 \,dx
 . \label{eq:thetaH1-3}
\end{align}
The second  integral in the r.h.s. of \eqref{eq:thetaH1-3} can be estimated  as
\begin{align*}
\int_{\bbS}\av{\tan(\theta) } \frac{(\pa_x\theta)^2\av{\pa_x^2\theta}}{1+\theta^2} \,dx&\le \eps \norm{\frac{\pa_x^2\theta(t)}{\sqrt{1+\theta^2(t)}}}_{L^2(\bbS)}^2+c(\eps)\norm{\pa_x\theta(t)}_{L^4(\bbS)}^4.
\end{align*}
Since
\[
\norm{\pa_x\theta(t)}_{L^4(\bbS)}^4\le \norm{\pa_x^2\theta(t)}_{L^2(\bbS)}^2\norm{\theta(t)}_{L^\infty(\bbS)}^2\le \norm{\frac{\pa_x^2\theta(t)}{\sqrt{1+\theta^2(t)}}}_{L^2(\bbS)}^2\norm{\theta(0)}_{L^\infty(\bbS)}^2\pare{1+\norm{\theta(0)}_{L^\infty(\bbS)}^2},
\]
we can absorb
\begin{align*}
\int_{\bbS}\av{\tan(\theta) } \frac{(\pa_x\theta)^2\av{\pa_x^2\theta}}{1+\theta^2} \,dx \le \pare{\pare{1+\norm{\theta(0)}_{L^\infty(\bbS)}^2}\norm{\theta(0)}_{L^\infty(\bbS)}^2+\eps} \norm{\frac{\pa_x^2\theta(t)}{\sqrt{1+\theta^2(t)}}}_{L^2(\bbS)}^2
\end{align*}
in the l.h.s. of \eqref{eq:thetaH1-3} for sufficiently small $\eps$ and  $\norm{\theta(0)}_{L^\infty(\bbS)}$.\\
The last integral in  the r.h.s. of \eqref{eq:thetaH1-3} can be estimated by Sobolev's embedding as
\[
\frac{1}{2}\int_{\bbS}\pare{
1+\frac{3+\theta^2}{1+\theta^2}\tan^2(\theta)
} \pare{\pa_{x}\theta}^3 \,dx  \le c\norm{\pa_x\theta(t)}_{L^3(\bbS)}^3\le  c\norm{\pa_x\theta(t)}_{H^\frac{1}{6}(\bbS)}^3.
\]
We interpolate the last inequality and we use the Poincaré inequality to get
\[
  c\norm{\pa_x\theta(t)}_{H^\frac{1}{6}(\bbS)}^3\le  c\norm{\pa_x\theta(t)}_{L^2(\bbS)}^\frac{5}{2}\norm{\pa_x^2\theta(t)}_{L^2(\bbS)}^\frac{1}{2}\le  c\norm{\pa_x\theta(t)}_{L^2(\bbS)}^2\norm{\pa_x^2\theta(t)}_{L^2(\bbS)}.
\]
We apply the Young inequality to obtain the following estimate on this last integral:
\[
\frac{1}{2}\int_{\bbS}\pare{
1+\frac{3+\theta^2}{1+\theta^2}\tan^2(\theta)
} \pare{\pa_{x}\theta}^3 \,dx \le c(\eps)\norm{\pa_x\theta(t)}_{L^2(\bbS)}^4+\eps\norm{\frac{\pa_x^2\theta(t)}{\sqrt{1+\theta^2(t)}}}_{L^2(\bbS)}^2.
\]
Finally, for $\eps$ and $\norm{\theta(0)}_{L^\infty(\bbS)}$ suitably small, i.e. such that
\[
\delta=1-3\eps-\pare{1+\norm{\theta(0)}_{L^\infty(\bbS)}^2}\norm{\theta(0)}_{L^\infty(\bbS)}^2>0,
\]
we obtain the following estimate:
\begin{align*}
 \frac{d}{dt}\int_{\bbS}u\pare{1+\tan^2 \pare{\theta}}\pare{\pa_{x}\theta}^2\,dx+
\delta\int_{\bbS}\frac{\pare{\pa_x^2\theta}^2}{1+\theta^2}\,dx  \le c \pare{\int_{\bbS}(\pa_x\theta)^2\,dx}\pare{\int_{\bbS}(\pa_x\theta)^2\,dx+1}.
\end{align*}
To conclude our argument, we continue in the following way: we estimate
\begin{align*}
 \int_{\bbS}(\pa_x\theta)^2\,dx &=\int_{\bbS}\frac{u\pare{1+\tan^2 \pare{\theta}}}{u\pare{1+\tan^2 \pare{\theta}}}(\pa_x\theta)^2\,dx
 \le c\int_{\bbS} u\pare{1+\tan^2 \pare{\theta}}(\pa_x\theta)^2\,dx.
\end{align*}
This means that we have a differential inequality of the type
\[
y'(t)\le c\pare{1+\pare{\int_{\bbS}(\pa_x\theta)^2\,dx}}y(t),
\]
with
\[
y(t)=\int_{\bbS}u\pare{1+\tan^2 \pare{\theta}}\pare{\pa_{x}\theta}^2\,dx.
\]
Then, we can apply a Grönwall type inequality obtaining that
\begin{align*}
&\int_{\bbS}u(t)\pare{1+\tan^2 \pare{\theta(t)}}\pare{\pa_{x}\theta(t)}^2\,dx\\
&\le \pare{\int_{\bbS}u_0\pare{1+\tan^2 \pare{\theta(0)}}\pare{\pa_{x}\theta(0)}^2\,dx}e^{c\pare{\int_0^t\pare{\int_{\bbS}(\pa_x\theta(s))^2\,dx+1}\,ds}}.
\end{align*}
Finally, we conclude using
$$
\int_0^t\int_{\bbS}(\partial_x \theta)^2\,dx\leq L(u_0).
$$

\section{Proof of Theorem \ref{teo:wiener} - Existence in Wiener spaces}
This proof is done in the same spirit as \cite{GMcristales,M,GMh,GMturra}. First we focus on the a priori estimates. We know that, for an $H^3$ initial data, the solution exists. Consequently, let us obtain the desired estimates under this extra assumption and, later on, we will generalize the argument to be able to drop it out.

Let
\[
w(x,t)=\frac{u(x,t)-\langle u_0\rangle}{\langle u_0\rangle}.
\]
Then, \eqref{eq:patu} becomes
\begin{equation*}%\label{eq:w}
\langle u_0\rangle w_t=\frac{\pa_x^2w+w\pa_x^2w-\pare{\pa_xw}^2}{1+2w+w^2+(\pa_xw)^2},
\end{equation*}
with initial data
\[
w_0(x)=w(x,0)=\frac{u_0(x)-\langle u_0\rangle}{\langle u_0\rangle}.
\]

Since we assumed that $\norm{w_0}_{(\bbS)}< 1/10$, we know that there exists a time $0<T^*$, eventually smaller than the existence time $T$, such that
\[
\norm{w(t)}_{A^1(\bbS)}< 1/10\qquad\forall t< T^*.
\]
Hence,
\[
\av{2w(x,t)+w^2(x,t)+(\pa_xw(x,t))^2 }\le 4\norm{w(t)}_{A^1(\bbS)}<1,
\]
and  we can develop in series
 \begin{equation}\label{eq:patw}
\langle u_0\rangle w_t=
\pa_x^2w+w\pa_x^2w-\pare{\pa_xw}^2+
\pare{\pa_x^2w+w\pa_x^2w-\pare{\pa_xw}^2}\pare{
\sum_{n\ge1}(-1)^n(2w+w^2+(\pa_xw)^2)^n
}.
\end{equation}
We want to write the Fourier coefficient of \eqref{eq:patw}. Before getting into these computations, we rewrite
\begin{align*}
\sum_{n\ge1}(-1)^n(2w+w^2+(\pa_xw)^2)^n&=\sum_{n\ge1}(-1)^n
\sum_{\substack{\ell+j+r=n\\ \ell,\, j,\, r\ge0}}
\binom{n}{\ell,\, j,\, r}2^\ell w^{2j+\ell}(\pa_xw)^{2r}
%\sum_{n\ge1}(-1)^n\sum_{j=0}^{n}\binom{n}{j}(\pa_xw)^{2(n-j)}\sum_{i=0}^{j}\binom{j}{i}2^iw^{2j-i},
\end{align*}
for
\[
\binom{n}{\ell,\, j,\, r}=\frac{n!}{\ell!\,j!\,r!}.
\]
The Fourier coefficient of this series is given by
\begin{align*}
&\w{\pare{\pa_x^2w+w\pa_x^2w-\pare{\pa_xw}^2}\pare{
\sum_{n\ge1}(-1)^n(2w+w^2+(\pa_xw)^2)^n
}}(k,t)\\
&= \w{\pare{\pa_x^2w+w\pa_x^2w-\pare{\pa_xw}^2}\pare{\ds\sum_{n\ge1}(-1)^n
\sum_{\substack{\ell+j+r=n\\ \ell,\, j,\, r\ge0}}
\binom{n}{\ell,\, j,\, r}2^\ell w^{2j+\ell}(\pa_xw)^{2r} }}(k,t)\\
&=\sum_{n\ge1}(-1)^n
\sum_{\substack{\ell+j+r=n\\ \ell,\, j,\, r\ge0}}
\binom{n}{\ell,\, j,\, r}2^\ell \pare{\w{w^{2j+\ell}(\pa_xw)^{2r}\pa_x^2 w}+\w{w^{2j+\ell+1}(\pa_xw)^{2r}\pa_x^2 w}-\w{w^{2j+\ell}(\pa_xw)^{2r+2}}} (k,t)
\\
&=\w{\m{N}}_1(k,t)+\w{\m{N}}_2(k,t)-\w{\m{N}}_3(k,t),
\end{align*}
for
\begin{align*}
\w{\m{N}}_1(k,t)&=-\sum_{n\ge1}(-1)^n
\sum_{\substack{\ell+j+r=n\\ \ell,\, j,\, r\ge0}}
\binom{n}{\ell,\, j,\,r}(-1)^{r}2^\ell\sum_{m_0\in\bbZ}\ldots\sum_{m_{2j+\ell+2r-3}\in\bbZ}(k-m_0)^2\w{w}(k-m_0,t)\times\\
&\qq\times (m_0-m_1)\w{w}(m_0-m_1,t)\pare{\prod_{p=1}^{2r-1}\w{w}(m_p-m_{p+1},t)(m_p-m_{p+1})}\times\\
&\qq \times \pare{\prod_{q=2r-1}^{2j+\ell+2r-3}\w{w}(m_q-m_{q+1},t)}\w{w}(m_{2j+\ell+2r-3},t),
\\
\w{\m{N}}_2(k,t)&=-\sum_{n\ge1}(-1)^n
\sum_{\substack{\ell+j+r=n\\ \ell,\, j,\, r\ge0}}
\binom{n}{\ell,\, j,\,r}(-1)^{r}2^\ell\sum_{m_0\in\bbZ}\ldots\sum_{m_{2j+\ell+2r-2}\in\bbZ}(k-m_0)^2\w{w}(k-m_0,t)\times\\
&\qq\times (m_0-m_1)\w{w}(m_0-m_1,t)\pare{\prod_{p=1}^{2r-1}\w{w}(m_p-m_{p+1},t)(m_p-m_{p+1})}\times\\
&\qq \times \pare{\prod_{q=2r-1}^{2j+\ell+2r-2}\w{w}(m_q-m_{q+1},t)}\w{w}(m_{2j+\ell+2r-2},t),
\\
\w{\m{N}}_3(k,t)&=\sum_{n\ge1}(-1)^n
\sum_{\substack{\ell+j+r=n\\ \ell,\, j,\, r\ge0}}
\binom{n}{\ell,\, j,\,r}(-1)^{r+1}2^\ell\sum_{m_0\in\bbZ}\ldots\sum_{m_{2j+\ell+2r-2}\in\bbZ}\times\\
&\qq\times (k-m_0)\w{w}(k-m_0,t)\pare{\prod_{p=0}^{2r}\w{w}(m_p-m_{p+1},t)(m_p-m_{p+1})}\times\\
&\qq \times \pare{\prod_{q=2r-1}^{2j+\ell+2r-2}\w{w}(m_q-m_{q+1},t)}\w{w}(m_{2j+\ell+2r-2},t),
\end{align*}
and the Fourier coefficient of \eqref{eq:patw} is given by
\begin{align*}
\langle u_0\rangle \w{w}_t(k,t)&=-k^2 \w{w}(k,t)-\sum_{m\in\bbZ}m^2\w{w}(m,t)\w{w}(k-m,t)
-\sum_{m\in\bbZ}m(k-m)\w{w}(m,t)\w{w}(k-m,t)\\
&\q +
\m{N}_1(k,t)+\m{N}_2(k,t)-\m{N}_3(k,t).
\end{align*}
We now want to take the sum in $k\in\bbZ$ and to estimate the $A^1(\bbS)$ semi-norm of the time derivative.\\
Since
\[
\pat|\w{w}(t,k)|={Re\left( \overline{\w{w}}(t,k)\pat \w{w}(t,k) \right)}/{|\w{w}(t,k)|},
\]
we have that
\begin{align*}
\sum_{k\in\bbZ}|k|\pat|\w{w}(t,k)|=\frac{d}{dt}\norm{w(t)}_{A^1(\bbS)}\qq\t{and}\qq
\sum_{k\in\bbZ}|k|^3\av{\w{w}(k,t)}=\norm{w(t)}_{A^3(\bbS)}.
\end{align*}
Using the Tonelli's Theorem and  interpolation in Wiener spaces, we estimate the left terms as
\begin{align*}
\sum_{k\in\bbZ}|k|\av{\sum_{m\in\bbZ}m^2\w{w}(m,t)\w{w}(k-m,t)}&\le \norm{w(t)}_{A^1(\bbS)}\norm{w(t)}_{A^2(\bbS)}\le  \norm{w(t)}_{A^0(\bbS)}\norm{w(t)}_{A^3(\bbS)},\\
\sum_{k\in\bbZ}|k|\av{\sum_{m\in\bbZ}m(k-m)\w{w}(m,t)\w{w}(k-m,t)}&\le \norm{w(t)}_{A^1(\bbS)}\norm{w(t)}_{A^2(\bbS)}\le \norm{w(t)}_{A^0(\bbS)}\norm{w(t)}_{A^3(\bbS)},
\end{align*}
and
\begin{align*}
\sum_{k\in\bbZ}|k|\av{\w{\m{N}}_1(k,t) }&\le \norm{w(t)}_{A^3(\bbS)}\sum_{n\ge1}\sum_{\substack{\ell+j+r=n\\ \ell,\, j,\, r\ge0}}\binom{n}{\ell,\, j,\,r}2^\ell\norm{w(t)}_{A^1(\bbS)}^{2r}\norm{w(t)}_{A^0(\bbS)}^{2j+\ell}\\
&=\norm{w(t)}_{A^3(\bbS)} \sum_{n\ge1}\pare{2\norm{w(t)}_{A^0(\bbS)}+\norm{w(t)}_{A^0(\bbS)}^2+\norm{w(t)}_{A^1(\bbS)}^2}^n,\\
\sum_{k\in\bbZ}|k|\av{\w{\m{N}}_2(k,t) }&\le \norm{w(t)}_{A^3(\bbS)}\sum_{n\ge1}\sum_{\substack{\ell+j+r=n\\ \ell,\, j,\, r\ge0}}\binom{n}{\ell,\, j,\,r}2^\ell\norm{w(t)}_{A^1(\bbS)}^{2r}\norm{w(t)}_{A^0(\bbS)}^{2j+\ell+1}\\
&=\norm{w(t)}_{A^3(\bbS)}\norm{w(t)}_{A^0(\bbS)} \sum_{n\ge1}\pare{2\norm{w(t)}_{A^0(\bbS)}+\norm{w(t)}_{A^0(\bbS)}^2+\norm{w(t)}_{A^1(\bbS)}^2}^n,\\
\sum_{k\in\bbZ}|k|\av{\w{\m{N}}_3(k,t) }&\le\norm{w(t)}_{A^2(\bbS)} \norm{w(t)}_{A^1(\bbS)} \sum_{n\ge1}\sum_{\substack{\ell+j+r=n\\ \ell,\, j,\, r\ge0}}\binom{n}{\ell,\, j,\,r}2^\ell\norm{w(t)}_{A^1(\bbS)}^{2r}\norm{w(t)}_{A^0(\bbS)}^{2j+\ell}\\
&\le \norm{w(t)}_{A^3(\bbS)} \norm{w(t)}_{A^0(\bbS)}\sum_{n\ge1}\sum_{\substack{\ell+j+r=n\\ \ell,\, j,\, r\ge0}}\binom{n}{\ell,\, j,\,r}2^\ell\norm{w(t)}_{A^1(\bbS)}^{2r}\norm{w(t)}_{A^0(\bbS)}^{2j+\ell}\\
&=\norm{w(t)}_{A^3(\bbS)}\norm{w(t)}_{A^0(\bbS)} \sum_{n\ge1}\pare{2\norm{w(t)}_{A^0(\bbS)}+\norm{w(t)}_{A^0(\bbS)}^2+\norm{w(t)}_{A^1(\bbS)}^2}^n.
\end{align*}
Since, for every $t<T^*$,
\begin{align*}
&\sum_{n\ge1}\pare{2\norm{w(t)}_{A^0(\bbS)}+\norm{w(t)}_{A^0(\bbS)}^2+\norm{w(t)}_{A^1(\bbS)}^2}^n\\
&=\frac{1}{1-(2\norm{w(t)}_{A^0(\bbS)}+\norm{w(t)}_{A^0(\bbS)}^2+\norm{w(t)}_{A^1(\bbS)}^2)}-1\\
&=\frac{2\norm{w(t)}_{A^0(\bbS)}+\norm{w(t)}_{A^0(\bbS)}^2+\norm{w(t)}_{A^1(\bbS)}^2}{1-(2\norm{w(t)}_{A^0(\bbS)}+\norm{w(t)}_{A^0(\bbS)}^2+\norm{w(t)}_{A^1(\bbS)}^2)}\\
&\le \frac{4\norm{w(t)}_{A^1(\bbS)}}{1-4\norm{w(t)}_{A^1(\bbS)}},
\end{align*}
we improve the bounds on the $\mathcal{N}_i$ terms as
\begin{align*}
&\sum_{k\in\bbZ}|k|\av{\w{\m{N}}_1(k,t) }+\sum_{k\in\bbZ}|k|\av{\w{\m{N}}_2(k,t) }+\sum_{k\in\bbZ}|k|\av{\w{\m{N}}_3(k,t) }\\
&\le  \norm{w(t)}_{A^3(\bbS)}\pare{1+2\norm{w(t)}_{A^0(\bbS)}} \frac{4\norm{w(t)}_{A^1(\bbS)}}{1-4\norm{w(t)}_{A^1(\bbS)}}.
\end{align*}
Then,
\begin{align*}
&\langle u_0\rangle \frac{d}{dt}\norm{w(t)}_{A^1(\bbS)}+\norm{w(t)}_{A^3(\bbS)}\\
&\le \pare{2\norm{w(t)}_{A^0(\bbS)}
+\pare{2\norm{w(t)}_{A^0(\bbS)}+1}\frac{4\norm{w(t)}_{A^1(\bbS)}}{1-4\norm{w(t)}_{A^1(\bbS)}}
}\norm{w(t)}_{A^3(\bbS)}\\
&= \pare{
\frac{2\norm{w(t)}_{A^0(\bbS)}}{1-4\norm{w(t)}_{A^1(\bbS)}^2}+\frac{4\norm{w(t)}_{A^1(\bbS)}}{1-4\norm{w(t)}_{A^1(\bbS)}}
}\norm{w(t)}_{A^3(\bbS)}\\
&\le \frac{6\norm{w(t)}_{A^1(\bbS)}}{1-4\norm{w(t)}_{A^1(\bbS)}}\norm{w(t)}_{A^3(\bbS)}.
\end{align*}
The smallness assumptions on $\norm{w_0}_{A^0(\bbS)}$ and on the $t<T^*$ imply that
\[
\frac{6\norm{w(t)}_{A^1(\bbS)}}{1-4\norm{w(t)}_{A^1(\bbS)}}\le c<1,
\]
so we get that
\[
\langle u_0\rangle \frac{d}{dt}\norm{w(t)}_{A^1(\bbS)}+(1-c)\norm{w(t)}_{A^3(\bbS)}\le 0,
\]
and the thesis follows integrating in time.

Let us explain now how to get the extra assumption on the initial data. We consider the problem
\begin{align*}
\langle u_0\rangle w_t^N&=
\pa_x^2w^N+w^N\pa_x^2w^N-\pare{\pa_xw^N}^2\\
&\q+
\pare{\pa_x^2w^N+w^N\pa_x^2w^N-\pare{\pa_xw^N}^2}\pare{
\sum_{n=1}^{N}(-1)^n(2w^N+(w^N)^2+(\pa_xw^N)^2)^n
}.
\end{align*}
For each fixed $N$ we can obtain a local in time solution using Galerkin method. Repeating the previous estimates in the Wiener space $A^1$, we find that
\begin{align*}
w^N\overset{*}{\rightharpoonup}  w\quad &\text{in} \quad L^{\infty}(0,T;L^\infty),\\%\label{linflinf}
\partial^3 w^N \overset{*}{\rightharpoonup}  \partial^3 w\quad& \text{in} \quad \mathcal{M}(0,T;L^\infty).%\label{MW4}
\end{align*}
Using the ideas in \cite{GGS} we can further get the desired space
$$
L^{\infty}(0,T;A^1)\cap \mathcal{L^1}(0,T;A^3).
$$

\end{document}